\documentstyle[12pt,leqno,amssymb]{amsart}

\newtheorem{thm}{Theorem}[section]
\newtheorem{lemma}[thm]{Lemma}
\newtheorem{prop}[thm]{Proposition}

\theoremstyle{definition}
\newtheorem{dfn}[thm]{Definition}
\theoremstyle{remark}

\begin{document}
\newcommand{\tc}{{\mathfrak{T}}}
\newcommand{\pr}{\protect\ref}
\newcommand{\ct}{\cite}

\newcommand{\h}{{1 \over 2}}
\newcommand{\su}{\subseteq}
\newcommand{\r}{{\mathrm{rank}}}
\newcommand{\C}{{\mathcal C}}
\newcommand{\F}{{\mathcal F}}
\newcommand{\ep}{{\epsilon}}
\newcommand{\A}{{\mathcal A}}

\newcommand{\lm}{{\lambda}}

\newcommand{\lf}{\lfloor}
\newcommand{\rf}{\rfloor}

\newcommand{\R}{{\Bbb R}}
\newcommand{\B}{{\Bbb B}}
\newcommand{\E}{{{\Bbb R}^3}}
\newcommand{\G}{{\Bbb G}}
\newcommand{\Z}{{\Bbb Z}}
\newcommand{\tf}{{\Bbb O}}

\newcommand{\im}{{Imm(F,\E)}}
\newcommand{\hc}{{H_1(F,\Z/2)}}

\newcounter{numb}

\title[Formulae for order one invariants]
{Formulae for order one invariants \\ of immersions and embeddings of surfaces}
\author{Tahl Nowik}
\address{Department of Mathematics, Bar-Ilan University, 
Ramat-Gan 52900, Israel.}
\email{tahl@@math.biu.ac.il}
\date{May 4, 2003}
\thanks{Partially supported by the Minerva Foundation}

\begin{abstract}
The universal order 1 invariant 
$f^U$
of immersions of a closed orientable surface into $\E$,
whose existence has been established in \ct{o}, takes values in the group 
$\G_U = K \oplus \Z/2 \oplus \Z/2$
where $K$ is a countably generated free Abelian group.
The projections of $f^U$ to $K$ and to the first and second $\Z/2$ factors 
are denoted $f^K, M, Q$ respectively. An explicit formula for the value of 
$Q$ on any embedding has been given in \ct{a}.

In the present work we give an explicit formula for the value of $f^K$ on any immersion, 
and for the value of $M$ on any embedding. 
\end{abstract}

\maketitle

\section{introduction}\label{intro}
Finite order invariants of stable immersions of a closed orientable
surface into $\E$ have been defined in \ct{o}, where all order 1
invariants have been classified. In \ct{h} all higher order invariants have been classified, 
and it has been shown that they are all functions of order 1 invariants.
This brings the attention back to order 1 invariants, and to the problem of finding explicit 
formulae for them. In \ct{o}, the existence of
a universal order 1 invariant $f^U$ has been established, 
which takes values in a group $\G_U = K \oplus \Z/2 \oplus \Z/2$ where $K$ is a 
countably generated free Abelian group. 
The existence proof, however, gave no clue
for computing the invariant.
We will denote the 
projections of $f^U$ to $K$ and to the first and second $\Z/2$ factors 
of $\G_U$ by $f^K, M, Q$ respectively. (The geometric meaning of $M$ and $Q$ will be explained
in Section \pr{inv}.) 
In \ct{a}, an explicit formula has been
given for $Q(i\circ h) - Q(i)$ where $h:F \to F$ is a diffeomorphism such that
$i$ and $i\circ h$ are regularly homotopic, and for
$Q(e')-Q(e)$ where $e,e'$ are any two regularly homotopic embeddings. 

In the present work we give an explicit formula for:
\begin{enumerate}
\item The value of $f^K$ on all immersions. 
\item $M(i\circ h) - M(i)$ where $h:F \to F$ is a diffeomorphism such that 
$i$ and $i\circ h$ are regularly homotopic. 
\item $M(e')-M(e)$ for any two regularly homotopic embeddings.
\end{enumerate}
Note that the invariant $f^U$ is specified only up to an 
order 0 invariant, i.e. up to an
additive constant in each
regular homotopy class, and so the same is true for $f^K,M,Q$.
For $M$ and $Q$ we will not have a specific choice of constants, 
and so as in (2),(3) above, we will speak only of the difference of the value
of $M$ and $Q$ on regularly homotopic immersions. 

The structure of the paper is as follows:
In Section \pr{back} we give the necessary background.
Note that in the present work we deviate from \ct{o},\ct{h} in our procedure for defining
order one invariants, and accordingly we deviate in our choice of generators for $\G_U$.
This is of no consequence in the abstract setting of \ct{o},\ct{h}, but will
greatly effect the simplicity of the explicit formula
for $f^K$ that we will find in the present work.
In Section \pr{inv} we explain the geometric meaning of the invariants
$M$ and $Q$. 
In Section \pr{st} we present the formulae that will be proved in this paper.
In Section \pr{fk} we prove the formula for $f^K$.
In Section \pr{ap} we give two applications. 
In Section \pr{m} we prove the formula for $M$.

\section{Background}\label{back}

In this section we summarize the background needed for this work.
Given a closed oriented surface $F$, $\im$ denotes the 
space of all immersions of $F$ into $\E$, with the $C^1$ topology.
A CE point of an immersion $i:F\to \E$ is a point of self intersection
of $i$ for which the local stratum in $\im$ corresponding to the 
self intersection, has codimension one.
We distinguish twelve types of CEs which we name
$E^0, E^1, E^2, H^1, H^2, T^0, T^1, T^2, T^3, Q^2, Q^3, Q^4$. 
Their precise description appears in the proof of Proposition
\pr{p1} below.
This set of twelve symbols is denoted $\C$.
A co-orientation for a CE is a choice of one of the two sides 
of the local stratum corresponding to the CE.
All but two of the above CE types are non-symmetric in the
sense that the two sides of the local stratum may be distinguished via
the local configuration of the CE, and for those ten CE types,
permanent co-orientations for the corresponding strata are chosen once and for all. 
The two exceptions are $H^1$ and $Q^2$ which are completely symmetric.
In fact, there does not exist a consistent choice of co-orientation for $H^1$ and $Q^2$
CEs since the global strata corresponding to these CE types are one sided in $\im$
(see \ct{o}).

We fix a closed oriented surface $F$ and a regular homotopy class $\A$ of 
immersions of $F$ into $\E$ (that is, $\A$ is a connected component of $\im$).
We denote by $I_n\su \A$ ($n\geq 0$) the space of all immersions in $\A$ which have precisely 
$n$ CE points (the self intersection being elsewhere stable).
In particular, $I_0$ is the space of all stable immersions in $\A$. 

Given an immersion $i\in I_n$, a \emph{temporary co-orientation} for $i$ is a choice 
of co-orientation at each of the $n$ CE points $p_1, \dots , p_n$ of $i$.
Given a temporary co-orientation $\tc$ for $i$ and a subset $A\su \{p_1,\dots,p_n\}$,
we define $i_{\tc,A} \in I_0$ to be the immersion obtained from $i$ by resolving all CEs
of $i$ at points of $A$ into the 
positive side with respect to $\tc$,
and all CEs not in $A$ into the negative side.
Now let $\G$ be any Abelian group and let $f:I_0\to\G$ be an invariant, i.e. a function 
which is 
constant on each connected component of $I_0$.
Given $i\in I_n$ and a temporary co-orientation $\tc$ for $i$,
$f^\tc(i)$ is defined as follows:
$$f^\tc(i)=\sum_{ A \su \{p_1,\dots,p_n\} } (-1)^{n-|A|} f(i_{\tc,A})$$
where $|A|$ is the number of elements in $A$.
The statement $f^\tc(i)=0$ is independent of the temporary co-orientation $\tc$
so we
simply write $f(i)=0$.
An invariant $f:I_0\to\G$ is called \emph{of finite order} if 
there is an $n$ such that $f(i)=0$ for all $i\in I_{n+1}$.
The minimal such $n$ is called the \emph{order} of $f$.
The group of all invariants on $I_0$ of order at most $n$ is denoted $V_n$.

From now on our discussion will reduce to order 1 invariants only.
The more general setting may be found in \ct{o},\ct{h}.
For an immersion $i:F\to\E$ and any $p\in\E$, we define the degree 
$d_p(i) \in \Z$ of $i$ at $p$ as follows: If $p \not\in i(F)$ then 
$d_p(i)$ is the (usual) degree of the map obtained from $i$ by composing it with the 
projection onto a small sphere centered at $p$. If on the other hand 
$p \in i(F)$ then we first 
push each sheet of $F$ which passes through $p$, a bit
into its preferred side determined by the orientation of $F$, obtaining
a new immersion $i'$ which misses $p$, and we define 
$d_p(i)=d_p(i')$.
If $i\in I_1$ and the unique CE of $i$ is located at $p \in \E$, then we define 
$C(i)$ to be the expression $R^a_m$ where $R^a\in\C$
is the symbol describing the configuration of the CE of $i$ at $p$ (one of the twelve symbols
above) and $m=d_p(i)$. We denote by $\C_1$ the set of all
expressions $R^a_m$ with $R^a\in\C, m\in\Z$. The map $C:I_1 \to \C_1$  
is surjective.

Let $f \in V_1$. For $i\in I_1$, if the CE of $i$ is
of type $H^1$ or $Q^2$ and $\tc$ is a temporary co-orientation for $i$,
then $2f^\tc(i)=0$ (\ct{o} Proposition 3.5),
and so in this case $f^\tc(i)$ is independent of $\tc$. 
This fact is used to extend any $f\in V_1$ to $I_1$
by setting for any $i\in I_1$, $f(i) = f^\tc(i)$,
where if the CE of $i$ is of type $H^1$ or $Q^2$ then
$\tc$ is arbitrary, and if it is not of type $H^1$ or $Q^2$
then the permanent co-orientation is used for the CE of $i$.
We will always assume without mention that any $f \in V_1$ is extended to $I_1$ 
in this way.
For $f\in V_1$ and $i,j\in I_1$, 
if $C(i)=C(j)$ then $f(i)=f(j)$ (\ct{o} Proposition 3.7),
so any $f\in V_1$ induces a well defined
function $u(f):\C_1\to\G$. 
The map $f\mapsto u(f)$ induces an injection $u:V_1 / V_0 \to \C_1^*$ 
where $\C_1^*$ is the group of all functions from $\C_1$ to $\G$.

The main result of \ct{o} is that the image of $u:V_1 / V_0 \to \C_1^*$ 
is the subgroup $\Delta_1 = \Delta_1(\G) \su \C_1^*$ which 
is defined as the set of functions in $\C_1^*$ satisfying 
relations which
we write as relations on the symbols $R^a_m$, e.g.
$T^0_m = T^3_m$ will stand for 
$g(T^0_m) = g(T^3_m)$.
The relations defining $\Delta_1$
are:
\begin{itemize}
\item $E^2_m = - E^0_m = H^2_m$, \ \  $E^1_m = H^1_m$.
\item $T^0_m = T^3_m$, \ \  $T^1_m = T^2_m$.
\item $2H^1_m =0 $, \ \ $H^1_m = H^1_{m-1}$.
\item $2Q^2_m =0 $, \ \ $Q^2_m = Q^2_{m-1}$.
\item $H^2_m - H^2_{m-1} = T^3_m - T^2_m$.
\item $Q^4_m - Q^3_m = T^3_m - T^3_{m-1}$, \ \  $Q^3_m - Q^2_m = T^2_m - T^2_{m-1}$. 
\end{itemize}

Let $\B\su\G$ be the subgroup defined by
$\B=\{ x\in \G : 2x=0\}$.
To obtain a function $g\in\Delta_1$ 
one may assign arbitrary values in $\G$ for the symbols
$\{T^2_m\}_{m\in\Z}$, $\{H^2_m\}_{m\in\Z}$ 
(here is where we deviate from \ct{o},\ct{h})
and arbitrary values in $\B$ for the two symbols $H^1_0 , Q^2_0$. 
Once this is done then 
the value of $g$ on all other symbols is uniquely determined, namely: 
\begin{enumerate}
\item $E^1_m = H^1_m = H^1_0$ for all $m$.
\item $E^2_m = -E^0_m = H^2_m$ for all $m$.
\item $T^3_m = T^2_m + H^2_m - H^2_{m-1}$
\item $T^0_m = T^3_m$, \ \ $T^1_m = T^2_m$ for all $m$.
\item $Q^2_m = Q^2_0$ for all $m$.
\item $Q^3_m (= Q^2_m + T^2_m - T^2_{m-1}) =Q^0_m + T^2_m - T^2_{m-1}$ for all $m$.
\item $Q^4_m (= Q^3_m + T^3_m - T^3_{m-1}) 
= Q^0_m + 2T^2_m - 2T^2_{m-1} + H^2_m - 2H^2_{m-1} + H^2_{m-2}$ 
for all $m$.
\end{enumerate}
In the sequel we will refer to this procedure as the "7-step procedure".

The Abelian group $\G_U$ is defined as follows (again note the difference from \ct{o},\ct{h}):
$$\G_U = \left< \{t^2_m\}_{m\in\Z}, \{h^2_m\}_{m\in\Z}, h^1_0, q^2_0 \ | \ 
2h^1_0 = 2q^2_0 = 0 \right>.$$ 

The universal element $g^U\in\Delta_1(\G_U)$ is defined by 
$g^U(T^2_m) = t^2_m$,
$g^U(H^2_m)=h^2_m$, $g^U(H^1_0)=h^1_0$, $g^U(Q^2_0)=q^2_0$ 
and the value of $g^U$ on all other symbols of $\C_1$ is determined by
the 7-step procedure. In \ct{o} the existence of
an order 1 invariant $f^U:I_0\to\G_U$ with $u(f^U)=g^U$ is proven. 
(Note that this is the same $g^U$ as in \ct{o} only presented via different generators).
The invariant $f^U$ is a \emph{universal} order 1 invariant, meaning the following:

\begin{dfn}\label{uni}
A pair $(\G,f)$ where $\G$ is an Abelian group and 
$f:I_0 \to \G$ is an order $n$ invariant,
will be called a \emph{universal order $n$ invariant}
if for any Abelian group $\G'$ and any order $n$ invariant 
$f':I_0 \to \G'$ there exists a unique homomorphism $\varphi:\G \to \G'$ 
such that $f' - \varphi \circ f$ is an invariant of order 
at most $n-1$. 
\end{dfn}

In \ct{h} all higher order invariants are classified, and for every $n$ a universal
order $n$ invariant is constructed as $\F_n \circ f^U$ where $\F_n:\G_U \to M_n$
is an explicit function (not homomorphism) into a certain Abelian group $M_n$.

\section{The invariants}\label{inv}

In this section we introduce the three invariants $f^K,M,Q$ that interest us.
We define $K \su \G_U$ to be the subgroup generated 
by $\{t^2_m\}_{m\in{\Bbb Z}} \cup \{h^2_m\}_{m\in{\Bbb Z}}$ 
(this is the same as the subgroup $K_1$ in \ct{h})
and define
$f^K:I_0 \to K$ to be the projection of $f^U$ to $K$. 
Similarly we define $M : I_0 \to \Z/2$ 
(respectively $Q: I_0 \to \Z/2$)
to be the projection of $f^U$ to $\Z/2$
sending all generators of $\G_U$ to 0 except $h^1_0$ 
(respectively except $q^2_0$).
Then $f^U = f^K \oplus M \oplus Q$.
Note that $f^U$ is defined only up to an additive constant in each regular homotopy class,
and so the same is true for $f^K,M,Q$.

More in detail, the invariants $Q$ and $M$ are defined as follows:
$M:I_0 \to \Z/2$ is the order 1 invariant defined by $u(M)(H^1_0)=1$, $u(M)(Q^2_0)=0$
and  $u(M)(T^2_m)=u(M)(H^2_m)=0$ for all $m$.
By the 7-step procedure, this extends to $u(M)(H^1_m)=u(M)(E^1_m)=1$ for all $m$,
$u(M)(H^a_m)=u(M)(E^a_m)=0$ for $a\neq 1$ and any $m$
and $u(M)(T^a_m)=u(M)(Q^a_m)=0$ for all $a,m$. That is, if $i_+,i_- \in I_0$ are the two 
immersions obtained from $i \in I_1$ by resolving its CE, then
$M(i_+)-M(i_-) = 1 \in \Z/2$ iff the CE of $i$ is a "matching tangency"
i.e. tangency of two sheets of the surface
where the orientations of the two sheets match at time of tangency.
(Thus the name $M$ for this invariant).
And so for any $i,j \in I_0$, $M(j)-M(i) \in \Z/2$ is the number mod 2 
of matching tangencies ocurring in any regular homotopy between $i$ and $j$.

Similarly, $Q:I_0 \to \Z/2$ is the $\Z/2$ valued order 1 invariant satisfying 
$u(Q)(Q^2_0) = 1$, $u(Q)(H^1_0)=0$
and $u(Q)(T^2_m)=u(Q)(H^2_m)=0$ for all $m$.
By the 7-step procedure, we have $u(Q)(Q^a_m)=1$ for all
$a,m$ and $u(Q)(T^a_m)=U(Q)(E^a_m)=u(Q)(H^a_m)=0$ for all $a,m$.
That is, $Q$ is the invariant such that for any $i,j \in I_0$, $Q(j)-Q(i)\in\Z/2$ is
the number mod 2 of quadruple points occurring in any regular homotopy between 
$i$ and $j$. This invariant has been studied in \ct{q} and \ct{a}.
In \ct{a} an explicit formula has been given for $Q(i \circ h) - Q(i)$
for any diffeomorphism $h:F \to F$ such that 
$i$ and $i \circ h$ are regularly homotopic,
and for $Q(e')-Q(e)$ for any two regularly homotopic embeddings. 
In the present work we will
do the same for $M$,
leaving open the interesting problem of finding an explicit formula for
$Q$ and $M$ on \emph{all} immersions. 
For $f^K$ however, we will indeed give a formula for all immersions.

\section{Statement of results}\label{st}

Let $i \in I_0$. For every $m\in\Z$ let $U_m = U_m(i) = \{ p \in \E-i(F) \ : \ d_p(i)=m \}$.
This is an open set in $\E$ which may be empty, and may be
non-connected or unbounded, but in any case, the Euler characteristic
$\chi(U_m)$ is defined.
Denote by $N_m = N_m(i)$ the number of triple points $p\in\E$ of $i$ having $d_p(i)=m$.  
The following formula for $f^K : I_0 \to K \su \G_U$ will be proved in Section \pr{fk}:

$$f^K(i)=
\sum_{m\in\Z} \chi(U_m) \bigg(\sum_{ -\h < k < \lf{m \over 2}\rf + \h} h^2_{m-2k}\bigg)
+ \sum_{m\in\Z} \h N_m \bigg(  t^2_m - \sum_{ -\h < k < m-\h } h^2_k \bigg)
$$
where for $a \in\R$, $\lf a \rf$ denotes the greatest integer $\leq a$,
and for $a,b \in \R$ the sum $\sum_{a<k<b}$ means the following:
If $a<b$ then it is the sum over all integers $a<k<b$, if $a=b$ then the sum is 0, and
if $a>b$ then $\sum_{a<k<b} = - \sum_{b<k<a}$.

For $i,j \in I_0$ let $M(i,j)=M(j)-M(i)$.
The following two formulae for $M$ will be proved in Section \pr{m}:

For any diffeomorphism $h:F\to F$ such
that $i$ and $i \circ h$ are regularly homotopic: 
$$M(i,i\circ h)=\bigg(\r(h_*-Id)\bigg)\bmod{2}$$
where $h_*$ is the map induced by $h$ on $H_1(F,\Z/2)$. 

If $e:F\to\E$ is an embedding then $e(F)$ splits $\E$ into two pieces,
one compact and one non-compact, 
which we denote $D^0(e)$ and $D^1(e)$ respectively.
By restriction of range, $e$ induces maps 
$e^k : F \to D^k(e)$, $k=0,1$. 
Let $e^k_* : \hc \to H_1(D^k(e),\Z/2)$ be the map 
induced by $e^k$.
Then for two regularly homotopic embeddings $e,e':F \to \E$,
$M(e,e')$ is computed as follows:
\begin{enumerate}
\item Find a basis 
$a_1,\dots,a_n,b_1,\dots,b_n$ for $\hc$ such that 
$e^0_*(a_i)=0$, $e^1_*(b_i)=0$ 
and $a_i \cdot b_j = \delta_{ij}$ (where $a \cdot b$ denotes the intersection form 
in $\hc$).
\item Find a similar basis
$a'_1,\dots,a'_n,b'_1,\dots,b'_n$ using $e'$ in place of $e$.
\item Let $m$ be the dimension of the subspace of $\hc$ spanned by:
$$a'_1 - a_1 \ , \ \dots \ ,  \ a'_n - a_n  \ ,  \ b'_1 - b_1  \ ,   \ 
\dots \ ,  \ b'_n - b_n.$$
\end{enumerate}
Then $M(e,e') = m\bmod{2} \in \Z/2$.

\section{Proof of formula for $f^K$}\label{fk}

We define the group $\tf$ to be the free Abelian group with generators 
$\{x_n\}_{n\in \Z} \cup \{y_n\}_{n\in\Z}$. 
For $i\in I_0$ we define $k(i) \in \tf$ as follows
(the terms are defined in Section \pr{st} and the sums are always finite):
$$k(i)= \sum_{m\in \Z} \chi(U_m) x_m + \sum_{m \in \Z} \h N_m y_m .$$
Indeed this is an element of $\tf$ since 
as we shall see below, $N_m$ is always even. In the mean time say
$k$ attains values in the $\Bbb Q$ vector space with same basis.

\begin{prop}\label{p1}
The invariant $k$ is an order 1 invariant, with $u(k)$ given by:
\begin{itemize}
\item $u(k)(E^a_m) = u(k)(H^a_m) = x_{m+a-2} - x_{m-a}$ 
\item $u(k)(T^a_m) = x_{m+a-3} + x_{m-a} + y_m$ 
\item $u(k)(Q^a_m) = x_{m+a-4} - x_{m-a} + (a-2)y_m + (2-a)y_{m-1}$ 
\end{itemize}
\end{prop}

\begin{pf}
We use the explicit description of the CE types, as appearing in \ct{o}, where more
details may be found.
A model in 3-space for the different sheets involved in the self intersection
near the CE, is given. The CE is obtained at the origin when setting $\lm=0$. 
We will show that for any $i \in I_1$, if $i_+ \in I_0$ 
is the 
immersion on the positive side of $i$ with respect to the permanent co-orientation 
for the CE of $i$, (if such exists, otherwise an arbitrary side is chosen)
and $i_- \in I_0$ is the immersion on the other side,
then indeed $k(i_+) - k(i_-)$
depends on $C(i)$ as in the statement
of this proposition. By showing in particular,
that this change depends \emph{only} on $C(i)$, we show
that $k$ is indeed an invariant of order 1.

Model for $E^a_m$: \ $z=0$, \ $z=x^2+y^2+\lm$.
The positive side is that where $\lm < 0$, where there is a new 2-sphere in the image
of the immersion, which is made of two 2-cells, and bounds a 3-cell in $\E$. 
The superscript $a$ is then the number of 2-cells (0, 1 or 2)
whose prefered side determined by the 
orientation of the surface, is facing away from the 3-cell, (and $m$ is the
degree at the CE at time $\lm=0$). The degree of points in 
the new 3-cell is seen to be $m+a-2$, and its $\chi$  is 1, and so the
term $x_{m+a-2}$. The second change ocurring, is that the region just above
the plane $z=0$, has a 2-handle removed from it, so its $\chi$ is reduced by
1, and the degree in this region is seen to be $m-a$, and so the term $- x_{m-a}$.

Model for $H^a_m$: \  $z=0$,  \ $z=x^2-y^2+\lm$. 
The positive side for $H^2$ is that where both sheets have their preferred side
facing toward the region that is between them near the origin.
For $H^1$ a positive side is chosen arbitrarily.
By rotating the configuration if necessary, say the positive side is where $\lm <0$.
The superscript $a$ then denotes the number of sheets (1 or 2) 
whose preferred side is facing toward the region that
is between the two sheets near the origin, when $\lm <0$. The changes ocurring in the 
neighboring regions when passing from $\lm > 0$ to $\lm < 0$ are that a 1-handle is removed
from the region $X$ just above the $x$ axis, and a 1-handle is added to the region $Y$
just below the $y$ axis. The degree of $X$ is seen to be $m+a-2$ and 
since a 1-handle is removed, $\chi(X)$ increases by 1 and thus the term $x_{m+a-2}$.
The degree of $Y$ is seen to be $m-a$, and since a 1-handle is added,
$\chi(Y)$ decreases by 1 and thus the term $-x_{m-a}$.

Model for $T^a_m$:  \ $z=0$, \ $y=0$,  \ $z=y+x^2+\lm$.
The positive side for this configuration is when $\lm < 0$, 
where there is a new 2-sphere in the image
of the immersion, which is made of three 2-cells, and bounds a 3-cell in $\E$. 
The superscript $a$ is the number of 2-cells (0, 1, 2 or 3)
whose prefered side is facing away from the 3-cell. The degree in the new 3-cell
is $m+a-3$ and its $\chi$  is 1 and so the term $x_{m+a-3}$.
The second change ocurring is that a 1-handle is removed from the region near the
$x$ axis having negative $y$ values and positive $z$ values. The degree of this region
is $m-a$ and since a 1-handle is removed, $\chi$ is increased by 1 and so the term 
$x_{m-a}$. The last change that effects the value of $k$ is that two triple points
are added, each of degree $m$. This increases $\h N_m$ by 1 and so the term $y_m$.

Model for $Q^a_m$:  \ $z=0$, \ $y=0$, \ $x=0$, \ $z=x+y+\lm$. 
On both the positive and negative side there is a simplex created near the origin,
and the positive side is that where the majority of the four sheets are facing away from the
simplex (and for $Q^2$ a positive
side is chosen arbitrarily). The superscript $a$ denotes the number of sheets
(2,3 or 4) facing away from the simplex created on the positive side, its degree thus
seen to be $m+a-4$. 
The simplex on the negative side has $4-a$ sheets facing away from it and so
its degree is $m-a$. So when moving from the negative to the positive side, a 3-cell ($\chi=1$)
of degree $m-a$ is removed and a 3-cell of degree $m+a-4$ is added, and so the terms
$x_{m+a-4} - x_{m-a}$. In addition to that, the degree of the four triple points of the
simplex changes. On the positive side there are $a$ triple points with degree $m$, (namely,
the triple points which are opposite the faces which are facing away from the simplex), 
and $4-a$ triple points with degree $m-1$. On the negative side the situation is reversed,
i.e. there are $4-a$ triple points with degree $m$ and $a$ triple points with degree 
$m-1$. So the total change in $N_m$ is $a-(4-a) = 2a-4$ and the total change in $N_{m-1}$ is
$(4-a)-a = 4-2a$ and so the terms $(a-2)y_m + (2-a)y_{m-1}$.

\end{pf}

We can now verify that indeed the values of $k$ are in $\tf$
i.e. no half integer coefficients appear (which means $N_m$ is always even). 
From Proposition \pr{p1} we see that the change in the value of $k$ is in $\tf$ 
along any regular homotopy, and so it is enough to
show that the value is in $\tf$ for one immersion in any given
regular homotopy class. Indeed, we show a bit more:

\begin{lemma}\label{l1}
Let $g$ be the genus of $F$. Any immersion $i:F\to\E$ is regularly homotopic
to an immersion $j$ with $k(j) = (2-g)x_0 + (1-g)x_{-1}$.
\end{lemma}

\begin{pf}
By \ct{p}, any immersion $i:F\to \E$ is regularly homotopic to an immersion
whose image is of one of two standard forms, either a standard embedding, 
or an immersion obtained from a standard embedding by adding a ring to it. 
(For what we mean by a "ring" see \ct{a}.) For an embedding $e$, $k(e)$
is either $(2-g)x_0 + (1-g)x_{-1}$ or $(2-g)x_0 + (1-g)x_1$, depending on whether
the preferred side of $e(F)$, 
determined by the orientation of $F$, is facing the compact or the non-compact side 
of $e(F)$ in $\E$, respectively. Now take an orientation reversing diffeomorphism
$h:F\to F$ such that $e\circ h$ is regularly homotopic to $e$, to see that both
values are attained. (Such $h$ exists by \ct{p}, take e.g. an $h$ that
induces the identity on $H_1(F, \Z/2)$.)
Now, a ring added to such embedding bounds a solid torus, whose $\chi$ is 0,
and the topological type and degree of the other two components remains the same, and so 
by the same argument as for an embedding, 
the two values are attained in this case too.
\end{pf}

We define a 
homomorphism $\varphi:\G_U \to \tf$ on generators as follows:

\begin{itemize}
\item $\varphi(h^2_m) = x_m - x_{m-2}$ 
\item $\varphi(t^2_m) = x_{m-1} + x_{m-2} + y_m$
\item $\varphi(h^1_0) = \varphi(q^2_0) = 0$
\end{itemize}

By Proposition \pr{p1}, $u(k) = u(\varphi \circ f^U)$ and so
$k=\varphi \circ f^U + c$ where $c \in \tf$ is a constant.
We define the following homomorphism $F:\tf \to K$
satisfying that $F \circ \varphi$ is the projection of $\G_U$ onto $K$,
and so $F \circ k = F \circ \varphi \circ f^U + F(c) = f^K + F(c)$.
By redefining $f^U$ as $f^U + F(c)$ we have $F \circ k = f^K$.
We define $F$ on generators of $\tf$ as follows (the notation
involved is defined in Section \pr{st}):
$$
F(x_m)  = \sum_{ -\h < k < \lf{m \over 2}\rf + \h} h^2_{m-2k} \ \ \ \ \ \ \ \ \ \ \ \ \ \ 
F(y_m)  = t^2_m - \sum_{ -\h < k < m-\h } h^2_k
$$
One checks directly that indeed $F \circ \varphi$
maps each generator of $K$ to itself. 
Since $\varphi$ is not surjective, there was a certain choice in the construction of $F$.
Indeed the image of $\varphi$
is the subgroup of $\tf$ of all elements
$\sum A_m x_m + \sum B_m y_m$ with $A_m,B_m \in \Z$ satisfying
$\sum_m A_{2m} = \sum_m A_{2m+1} = \sum_m  B_m$. And so any two generators
$x_i , x_j$ with $i$ even and $j$ odd, generate a subgroup 
in $\tf$ which is a direct summand of the
image of $\varphi$. Our choice for $F$ was that $F(x_{-2}) = F(x_{-1}) = 0$.
Note that by Lemma \pr{l1}, the image of $k:I_0 \to \tf$ 
is contained in a non trivial coset of the image of $\varphi$ in $\tf$,
(and so the constant $c$ appearing above is non-zero, regardless of an additive 
constant for $f^U$).
Composing the formula for $F$ with the formula for $k$ we obtain our formula for $f^K$: 
$$f^K(i)=
\sum_{m\in\Z} \chi(U_m) \bigg(\sum_{ -\h < k < \lf{m \over 2}\rf + \h} h^2_{m-2k}\bigg)
+ \sum_{m\in\Z} \h N_m \bigg(  t^2_m - \sum_{ -\h < k < m-\h } h^2_k \bigg).
$$
The choice of constants for $f^K$ here may be characterized by saying that
in each regular homotopy class, $ f^K(j) = (2-g)h^2_0 $ for $j$ of Lemma \pr{l1}.

Since $f^U$ is universal, the image of $f^K:I_0\to K$ 
is not contained in any coset of a proper subgroup of $K$, yet
the image of $f^K$ is far from being the whole
group $K$, since as we see from the formula, the coefficients of all
generators $t^2_m$ are always non-negative. It would be interesting to
determine the precise image of $f^U:I_0 \to \G_U$.

\section{Applications}\label{ap}

We give the following two applications. The first will be used in the second
and the second will be used in Section \pr{m}.

We will use the fact that $\varphi:\G_U \to \tf$ is not surjective to obtain
identities on immersions:
Let $\theta_0,\theta_1 : \tf \to \Z$ be the homomorphisms defined by:
$\theta_0(x_{2m})=1, \theta_0(x_{2m+1})=0, \theta_0(y_m)=-1$ for all $m$ and 
$\theta_1(x_{2m})=0, \theta_1(x_{2m+1})=1, \theta_1(y_m)=-1$ for all $m$ and 
so $\theta_0 \circ \varphi = \theta_1 \circ \varphi = 0$. It follows that
$\theta_0 \circ k$ and
$\theta_1 \circ k$ are constant invariants, which are given explicitly by
$\theta_0 \circ k(i) = \sum \chi(U_{2m}) - {1 \over 2}N$
and $\theta_1 \circ k(i) = \sum \chi(U_{2m+1})  - {1 \over 2}N$
where $N=N(i)=\sum N_m(i)$ is the total number of triple points of $i$.
To find the value of these constants we need to evaluate them on a single immersion
in every regular homotopy class. For the immersion $j$ of Lemma \pr{l1}, 
$\theta_0 \circ k(j)=2-g$ and $\theta_1 \circ k(j)=1-g$,
so we get the following two identities:
For any $i \in I_0$, $$\sum_m \chi(U_{2m}) - {1 \over 2}N = 2-g \ \ \ \ \ \
\text{and} \ \ \ \ \ \  
\sum_m \chi(U_{2m+1})  - {1 \over 2}N = 1-g.$$
 
For our second application,
let $U:I_0 \to \Z$ be the order one invariant defined by 
$u(U)(H^2_m)=1$, $u(U)(T^2_m)=0$ for all $m$ and
$u(U)(H^1_0)=u(U)(Q^2_0)=0$.
By the 7-step procedure we will have 
$u(U)(H^2_m)=u(U)(E^2_m)=-u(U)(E^0_m)=1$ for all $m$,  
$u(U)(H^1_m)=u(U)(E^1_m)=0$ for all $m$,
and $u(U)(T^a_m)=u(U)(Q^a_m)=0$ for all $a,m$.
That is, for any $i,j \in I_0$,
$U(j)-U(i) \in \Z$ is the signed number of \emph{un}-matching tangencies occurring in
any regular homotopy from $i$ to $j$
(thus the name $U$ for this invariant) where each such tangency is counted as
$\pm 1$ according to its permanent co-orientation and the prescription
$u(U)(H^2_m)=u(U)(E^2_m)=-u(U)(E^0_m)=1$.
Following the definition of $U$ we define $\eta : K \to \Z$ on generators as follows:
$\eta(h^2_m)=1$ and $\eta(t^2_m)=0$ for all $m$.
Then $u(U) = u(\eta \circ f^K)$ and so (up to choice of constants) $U=\eta \circ f^K$. 
So from our formula for $f^k$ we get an explicit formula for $U$:
$$U(i)=\sum_{m\in\Z} \chi(U_m) \lf {m+2 \over 2} \rf - \sum_{m \in\Z} {1\over 2}m N_m .$$ 
Again we may characterize the choice of constants by saying that $U(j)=2-g$ for
$j$ of Lemma \pr{l1}

We denote $U(i,j)=U(j)-U(i)$.
For two regularly homotopic embeddings $e,e':F \to \E$ we would like
to compute $U(e,e')$.
For $e:F \to \E$ an embedding let 
$c(e)\in \Z$ be the degree of the points in the compact side of $i(F)$ in $\E$,
so $c(e)=\pm 1$.
We have 
$U(e) = (2-g) + (1-g) \lf {c(e) + 2 \over 2} \rf$ and so
$$U(e,e') = U(e')-U(e) 
= (1-g)\bigg(\lf {c(e') + 2 \over 2} \rf - \lf {c(e) + 2 \over 2} \rf\bigg)
= (1-g)\ep(e,e')$$
where $\ep(e,e')$ is
$0$ if $c(e)=c(e')$, is $1$ if $c(e)=-1,c(e')=1$ and is $-1$ 
if $c(e)=1,c(e')=-1$.

Now for $i \in I_0$ and $h:F \to F$ a diffeomorphism such
that $i$ and $i \circ h$ are regularly homotopic, we would like to compute 
$U(i,i\circ h)$. If $h$ is orientation preserving then from the formula we have for 
$U(i)$ it is clear that $U(i)=U(i \circ h)$ and so $U(i,i \circ h)=0$.
Now let $h:F \to F$ be orientation reversing. If $p \in \E - i(F)$ then
$d_p(i \circ h) = -d_p(i)$ and if $p\in\E$ is a triple point of $i$
then $d_p(i \circ h) = 3-d_p(i)$ and so we get:
$$U(i \circ h) - U(i)=
\sum_m \chi(U_m(i)) (\lf {-m+2 \over 2} \rf - \lf {m+2 \over 2} \rf) 
- \sum_m {1\over 2}(3-m - m) N_m(i) .$$ 
Using the two identities from the beginning of this section and the fact that 
$\lf {-m+2 \over 2} \rf - \lf {m+2 \over 2} \rf = -2\lf {m+2 \over 2} \rf + k(m)$ where
$k(m)$ is 2 for $m$ even and 1 for $m$ odd, we get:
$$U(i, i\circ h) = (1-g) + 2\bigg(2-g-U(i)\bigg).$$
Note that the $U(i)$ appearing here on the right, stands for our specific formula for
the invariant $U$, and not for the abstract invariant which is defined only up to a constant.
This equality for $h$ orientation reversing can be interpreted as 
$U(i,i\circ h) = U(j,j\circ h) + 2U(i,j)$ for $j$ of Lemma \pr{l1}, 
offering another way for proving the equality.

Let $\widehat{U}:I_0\to \Z/2$ be the mod 2 reduction of $U$. 
The reduction mod 2 of the above results reads as follows:
For embeddings $e,e':F \to \E$, $\widehat{U}(e,e')=(1-g)\widehat{\ep}(e,e')$
where $\widehat{\ep}(e,e')\in\Z/2$ is 0 if $c(e)=c(e')$ and is 1 if $c(e) \neq c(e')$. 
For $h:F\to F$ a diffeomorphism such that $i$ and $i\circ h$ are regularly homotopic,
$\widehat{U}(i,i\circ h)=(1-g)\ep(h)$ where $\ep(h) \in \Z/2$ is 0 if
$h$ is orientation preserving and is 1 if $h$ is orientation reversing.

\section{Proof of formula for $M$}\label{m}

For $i \in I_0$ and $h:F \to F$ a diffeomorphism such that $i$ and $i \circ h$ are
regularly homotopic, let $M'(i,i\circ h)$ 
denote our proposed formula for $M(i,i\circ h)$
presented in Section \pr{st}. So we must show that indeed 
$M(i,i\circ h)=M'(i,i\circ h)$.
Similarly, for regularly homotopic embeddings $e,e':F \to \E$, let $M'(e,e')$ denote
the proposed value for $M(e,e')$ presented in Section \pr{st}, so we must show
$M(e,e')=M'(e,e')$.

In \ct{a} it is shown that $Q(i,i\circ h) = M'(i,i\circ h) + (1-g)\ep(h)$ and
$Q(e,e')= M'(e,e') + (1-g)\widehat{\ep}(e,e')$.
In view of the concluding paragraph of Section \pr{ap}, this means
$Q(i,i\circ h) = M'(i,i\circ h) + \widehat{U}(i,i\circ h)$ and
$Q(e,e')= M'(e,e') + \widehat{U}(e,e')$.
So showing $M=M'$ in these two settings is equivalent to showing
$Q=M+\widehat{U}$ in these settings, which means that the number 
mod 2 of quadruple points occurring in any regular homotopy 
between such two immersions or embeddings, is equal to the 
number mod 2 of all tangencies occurring (matching and un-matching).
So, it remains to prove the following:

\begin{prop}\label{p2}
Let $i,j \in I_0$ such that either there is a diffeomorphism $h:F\to F$ such
that $j= i \circ h$ or $i,j$ are both embeddings. Then in any regular
homotopy between $i$ and $j$, the number mod 2 of quadruple points occurring,
is equal to the number mod 2 of tangencies occurring.
\end{prop}

\begin{pf}
For a closed 3-manifold $N$ and stable immersion $f:N \to \R^4$, there is defined a closed
surface $S_f$ and immersion $g:S_f \to \R^4$ such that the image $g(S_f) \su \R^4$
is precisely the multiple set of $f$. It is shown in \ct{ec} that the number mod 2
of quadruple points of $f$ is equal to $\chi(S_f) \bmod 2$.

Now let $i,j:F \to \E$ be as in the assumption of this proposition and
let $H_t:F \to \E$, $0 \leq t \leq 1$, be a regular homotopy with $H_0=i$, $H_1=j$.
We define an immersion $f:F \times [0,1] \to \E \times [0,1]$ by $f(x,t) = (H_t(x) , t)$.

In case $i,j$ are embeddings we continue $f$ into $\R^4 = \E \times \R$
and construct a closed 3-manifold $N$ by attaching
two handle bodies to $F \times [0,1]$, glued so that $f$ can be extended to embeddings
of these handle bodies into $\E \times (-\infty ,0]$ and $\E \times [1,\infty)$.
We thus obtain an immersion $\bar{f}:N\to\R^4$ with self intersection being 
precisely the original self intersection of $F \times [0,1]$. The projection 
$\E \times [0,1] \to [0,1]$ induces a Morse function on $S_{\bar{f}}$ with
singularities precisely wherever a tangency CE occurs in the regular homotopy $H_t$,
and so by Morse theory $\chi(S_{\bar{f}})$
is equal mod 2 to the number of tangencies.
By \ct{ec} then, the number mod 2 of quadruple points of $H_t$ which is the number mod 2
of quadruple points of $\bar{f}$ is equal to the number mod 2 of tangencies.

In case $j = i \circ h$, let $N$ be the 3-manifold obtained from $F \times [0,1]$ by
gluing its two boundary components to each other via $h$ so that there is induced an immersion
$\bar{f}:N \to \E \times S^1$. Composing $\bar{f}$ with an 
embedding of $\E \times S^1$ in $\R^4$, we see again that
the number of quadruple points of $H_t$ is equal mod 2 to $\chi(S_{\bar{f}})$
which is equal mod 2 to the number of tangencies of $H_t$.

\end{pf}

We remark that one can prove the formulae for $M$ presented in Section \pr{st}
directly, without resorting to the result of \ct{ec}, by going along the lines
of \ct{a}. Proposition \pr{p2} would then be obtained as a corollary.

\end{document}